\documentclass[11pt
]{article}

\usepackage[margin=0.75in]{geometry} 

\usepackage{titlesec}
\titleformat{\section}[runin]
{\normalfont\bfseries}
{\thesection.}{0.4em}{}

\usepackage{amsmath,amssymb,amsthm,latexsym}
\usepackage{graphicx}
\usepackage{float}
\usepackage{url}
\usepackage{combelow}
\usepackage{mathrsfs}
\usepackage[scr=rsfs,cal=boondox]{mathalfa}

\newtheorem*{theorem*}{Theorem}

\usepackage[
]{hyperref}

\usepackage{wrapfig}

\begin{document}

\title{A Wild Steiner-Lehmus Chase \\  {\small Being A Side-By-Side Proof of the Steiner-Lehmus Equal Bisectors Theorem} }
\markright{Notes}
\author{Eric L. Grinberg and Mehmet Z.Orhon}

\date{\vspace{-5ex}}

\maketitle

\begin{abstract}
We present a proof  the Steiner-Lehmus equal bisectors theorem by applying the Law of sines in rapid succession to a side-by-side comparison.  For nearly two centuries, the quest for a direct proof has sustained interest in proving and reproving this theorem.  We suggest that a second driving force may also be at play.

\end{abstract}

\section{ Every Proof Tells a Story}

 In the last days of 2021 we read Hillel Furstenberg's Abel prize interview \cite{dundas}  and watched the related video \cite{furstenberg}. He tells:
 \begin{quote}
\emph{\small There was a particular problem that was responsible for my career.}
 \end{quote}
 
 Early in his education, Furstenberg was challenged to prove that, in a triangle, if two of the angle bisectors have the same length then the triangle is isosceles. He succeeded, and that launched his career. Moreover, Shlomo Sternberg, Furstenberg classmate and future prominent mathematician,  had the same experience at the same time with the same problem; a later colleague of Furstenberg’s,  the mathematical economist and future Nobel laureate Robert J. Aumann had the same career launching experience as well. 
 
 Not long after seeing the Abel interviews, we watched a segment of the Simons Foundation’s \emph{Science Lives} series, this one featuring the distinguished mathematician Elias M.~Stein, interviewed by his former student Charles Fefferman \cite{fefferman-2012}. Minutes into the interview we learn that, early in his education,  Eli was  presented with the same theorem, along with an ``ingenious’’ proof.  He saw a   challenge and declared ``I’m gonna try to find another proof of this theorem.” After some time, success came:
 
 \vspace{-5pt}
 
 \begin{quote}
\emph{It instilled in me the confidence that  I could do mathematics.}
 \end{quote}
Could this be a coincidence? How many other careers did this theorem, and its (re)proofs launch?

We soon learned that the theorem has a name, the \emph{Steiner-Lehmus Theorem} (SLT for short), along with an extensive, interesting history, dotted with proofs and reproofs. Although the authors of the present paper are farther along in their careers than the aforementioned mathematical giants were, the challenge to re-prove appeared strong.  Awareness of another much-reproved theorem, \cite{oakley-1978}
%
%
 added activation energy. This led to a new elementary proof, yet another, which will be detailed below, after re-proof remarks.
 \\
 {\texttt { \footnotesize 
 [Added in re-proof:  As we were about to finalize paper, we learned that another prominent mathematician, \\ 
 Steven Strogatz, also experienced a Steiner-Lehmus career takeoff. See \cite{sanderson-2021}, where Grant Sanderson interviews Strogatz in the 3Blue1Brown Podcast series. Interestingly, Strogatz also brings up an \\ 
 inspiring proof he learned as an undergraduate, in a class taught by Eli Stein. This was Goursat's re-proof of the Cauchy Integral Theorem.]}

\section{Theorems That Keep On Giving.}
Soon after being introduced to the Steiner-Lehmus theorem we  discovered the Susan D’Agostino article  with the intriguing title \emph{Mathematicians Will Never Stop Proving the Prime Number Theorem} \cite{dagostino-2020}. We found particular interest in the subtitle:

\begin{center} 
\emph{Why do mathematicians enjoy  proving the same results in different ways?}
\end{center}

While the focus in \cite{dagostino-2020} is on the Prime Number Theorem, other, much re-proved theorems are also mentioned. We then  learned that there is an entire book about the re-prove practice: J.W.~Dawson’s 2015 work with the intriguing title \emph{Why Prove It Again?} \cite{dawson-2015}.  In addition, Paul Erd\H{o}s' celebrated concept of \emph{The Book}, famously articulated in \cite{aigner-2018}, is all about re-proofs, albeit only the best ones. In \cite{loomis-1968} we have an entire book about reproofs of a single theorem. Also relevant is D.~Hilbert's 24th problem and the related concept of \emph{purity of the method}, seeking the simplest possible proofs; see \cite{greenberg-2010}.

There is more to say about the paradigm, but we digress. It’s time for a (re)-proof.

\section{Proving The Steiner-Lehmus Theorem.}

Influenced by the excellent article \cite{bauer-2017},  we state the SLT in a form which eschews one contrapositive. 
 
\begin{theorem*}
Let the triangle $T \equiv \triangle ABC$ have base angles $\angle A, \angle B$, with $\angle A$ larger than $\angle B$. 
Then the internal bisecting segment of $\angle A$ in $T$ is shorter than the internal bisecting segment of $\angle B$ in $T$. 
\end{theorem*}
\noindent
In brief:
\qquad\qquad\qquad
\emph{The larger angle has the smaller bisector.} \\

\noindent\emph{Proof.} 
Begin with figure \ref{fig:figure01}, marking  angle measures $2 \alpha, 2\beta, \gamma$, corresponding to the vertices $A,B,C$, respectively. Bisecting the angle $\angle A$ we extend the bisecting segment until it meets the opposite side, $BC$, denoting both the resulting segment \emph{and} its length  by $\ell$. (We admit to abuse of notation here.)  Our base assumption is that $2\beta$ is the smaller of the two bottom angles:
\[
2\beta < 2\alpha \, .
\]

Now $\ell$ divides $\triangle ABC$ into two sub-triangles, which we will regard as \emph{$\ell$-upper} and \emph{$\ell$-lower}, respectively. In the $\ell$-lower triangle, the angle where $\ell$ meets $BC$ is supplementary to an angle on the $\ell$-upper triangle which is supplementary to $\gamma + \alpha$, so this angle equals $\gamma + \alpha$ (later to be denoted by $\gamma^\prime$):

\begin{figure} [H]
\caption{Triangle with first bisector}
\vspace{-90pt}
\centerline{\includegraphics[width=0.9\linewidth]{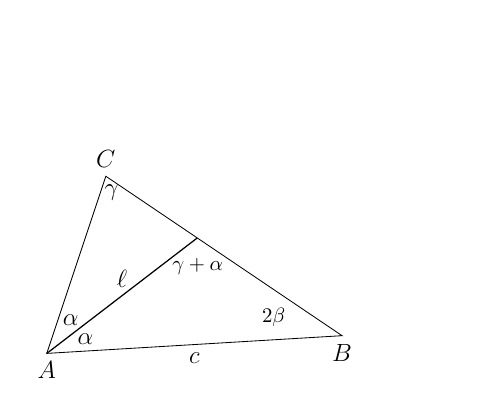}}
\label{fig:figure01}
\end{figure}
\noindent Admittedly, introducing the synonym $\gamma^{\prime}$ for $\gamma + \alpha$ adds notational baggage. But it also declutters ensuing diagrams.

Next, we bisect $\angle B$ and denote  both the corresponding bisecting segment and its length by $\ell^\prime$ (continuing the practice of mild abuse of notation). This  segment also divides $\triangle ABC$ into two sub-triangles, to be regarded as  \emph{$\ell^\prime$-upper} and \emph{$\ell^\prime$-lower}, respectively, as in figure \ref{fig:figure02}.

\begin{figure} [H]
\caption{Triangle with second bisector}
\centerline{\includegraphics[width=0.8\linewidth]{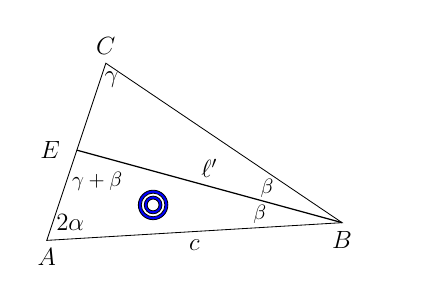}}
\label{fig:figure02}
\end{figure}

 Denote by $E$ the point where the bisecting segment $\ell^\prime$ meets the side $AC$. Adorn the lower $ \ell^\prime$ triangle ($\triangle ABE$) with a virtual blue and white ``bullseye" disk, which we will use to focus attention  when we apply  reflection and translation to that triangle. \\
%
%

\begin{wrapfigure}
[3]{r}{0.17\textwidth} 
\begin{center}
\vspace{-37pt} \!\!
\includegraphics[width=0.18\textwidth] {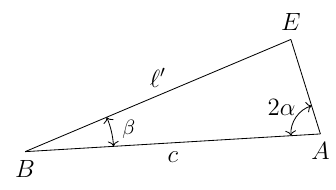}
\end{center}
\end{wrapfigure}

Take a copy of \emph{$\ell^\prime$-lower} triangle $\triangle AEB$ 
in figure \ref{fig:figure02}
and reflect it about the perpendicular to $AB$  through the bullseye,  obtaining a triangle similar to the miniature depicted on the right. 
Translate the reflected
triangle, placing it ``on top" of the original triangle $\triangle ABC$ of figure
\ref{fig:figure01},  so that the transported vertex $B$ is made to coincide with original vertex $A$, and the transported side $BA$ now lies along the original side $AB$, as in  figure \ref{fig:figure03}.

\begin{figure} [!ht]
\caption{Triangles, upon reflection}
\vspace*{-28mm}
\centerline{\includegraphics[width=0.9\linewidth]{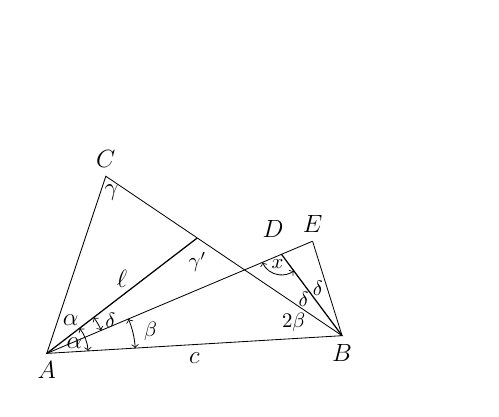}}
\label{fig:figure03}

\vspace{-20pt}
\qquad \qquad \qquad\qquad\qquad {
\flushright \hfill
\fbox{\begin{minipage}[b]{1.4in}
{  \smallskip 
\flushright 
\emph{I want you to compare \, }
\\ 
\bigskip
\hfill
\scriptsize  Lee Iacocca \\
Chrysler recovery campaign \\ \hfill  1980s \cite{iacocca}
  }
 \end{minipage}
 } }
\end{figure}

In figure \ref{fig:figure03} we use the notation  $\delta \equiv \alpha - \beta$, trading clutter for notational baggage oncemore.
 At vertex $B$,  mark an angle $2\beta + \delta$ from the base $AB$. This angle meets $AE$ at a point $D$. The angle $\angle x \equiv \angle ADB$ is supplementary to $\beta+2 \beta + \delta={\alpha+2 \beta}$. Now $\gamma^\prime$ is also  supplementary to ${\alpha+2 \beta}$, as the \emph{$\ell$-lower} triangle in figure \ref{fig:figure01} shows, so $\angle x = \gamma^\prime$.
 
 Applying the Law of sines to $\angle D \equiv \angle x \equiv \gamma^\prime$ and $\angle B$ in $\triangle ADB$ we have:
\[
\sin( \gamma^\prime )/c= \sin(2 \beta + \delta ) / AD \, .
\]
In the \emph{$\ell$-lower} triangle of figure \ref{fig:figure01}, the Law of sines applied to $\angle \gamma^\prime = \gamma+\alpha$ and $\angle 2 \beta$ gives
\[
\sin ( \gamma^\prime ) / c = \sin ( 2 \beta ) / \ell \, .
\]
Combining the last two equations we have
\begin{equation} \label{eq:first_estimate}
AD= \left(   \frac{ \sin(2 \beta + \delta )}{ \sin ( 2 \beta )} \right)  \ell 
%
%
\, .
\end{equation}

Now the authors of this paper may be  obtuse, but the angle $\angle (2 \beta + \delta)$ cannot be. This is because $\angle (2 \beta + \delta)$ is ``halfway" between $\angle 2 \beta$ and $\angle 2 \alpha$, two angles of a triangle whose remaining angle, $\gamma$ needs a little room to breathe.  Thus the sine function is strictly monotone increasing from $2 \beta$ to  $ 2 \beta + \delta $ and \eqref{eq:first_estimate} shows that $AD$ is longer than $\ell$. As $AD$ is a subsegment of $AE$, and $AE$ has length $\ell^\prime$, we infer that
\[
\ell^\prime > AD  > \ell
%
\, .
\]
So \emph{the larger angle has the smaller bisector} and this proves our statement of the  Steiner-Lehmus theorem. 
\qedsymbol 

Is the reflection step absolutely necessary? No, the proof may be re-couched without it, but reflection  does provide a \emph{visualization} tool--eyeball comparison of  the lengths of two segments is easier when the segments are more closely oriented. Indeed,  J.~Steiner's original proof also uses reversal as a visualization tool: ``\"Ubersicht'';  see \cite[p. 225]{ostermann-2012}.

We note that this proof is shorter than the Law of sines proof given by V.~Cristescu in 1916 \cite{cristescu-1916}, which is reproduced in \cite{olah-sandor-2009} and is there deemed to be among ``the shortest trigonometric proofs" of the SLT. Also, this proof is largely geometric with little in the way of algebra or analysis often included in trigonometric proofs.

\section{The Oldest Indirect Proof?}

Early in Euclid’s Elements we find the celebrated \emph{Pons Asinorum} theorem: in an isosceles triangle, the base angles are equal; the Euclid proof  is complicated, hence the name \cite[p. 105]{moise-1990}, but it is a direct proof. Today’s expositions typically invoke the later and considerably simpler proof by Pappus, which is also direct.  But Pons Asinorum has a converse: in a triangle, if two base angles are equal, then the corresponding sides are equal, i.e., the triangle is isosceles. This is proposition 6 in the Elements, and its proof is indirect from the start. Indeed it appears to be among the oldest popular indirect proofs. 

Modern expositions critique Euclid's proof of the converse to Pons Asinorum, as it implicitly relies on area considerations without developing a theory for such. In a modern rendition \cite[p. 98]{hartshorne-2000} the proof is almost the same as in the Elements, merely replacing the ending area argument with the Hilbert angle replication axiom, C4; it is still indirect.  Alternatively, 
\cite{moise-1990} takes a different approach, using the ASA congruence theorem, with indirectness integral; see below. Either way, it seems that the Steiner-Lehmus theorem has an ancient cousin with indirect, non-constructive roots.

\section{Constructive Criticism and Indirect Costs}

We have already declared that \cite{bauer-2017}  influenced us to present the main theorem in a form that avoids one use of contraposition in the proof. However, there are additional issues.

In preparing federal grant budgets there is the line item of \emph{indirect costs}. The same phrase applies here, though with a different meaning. We use tools from elementary geometry whose foundations  have indirect, non-constructive roots. One, which is needed throughout our proof, is that the total interior angle sum of a triangle is $180^{\circ}$ (or $\le 180^{\circ}$ in neutral geometry). In modern treatments this is proved indirectly, using analysis of \emph{Saccheri quadrilaterals}. See \cite[section 10.4]{moise-1990} Another is the celebrated \emph{ASA} triangle congruence theorem, which is proven non-constructively from the \emph{SAS} congruence principle. The SAS congruence, in turn, is a theorem in Euclid but an axiom in modern treatments; see \cite[chapters 6 and 10]{moise-1990}. Moreover, SAS and ASA figure prominently in the proof of consistency of Euclidean geometry. Thus, in surveying the multitude of proofs of the SLT, we find indirectness not only in basic structure, but also in reliance on indirectly sourced tools.

\section{D\`ej\'a Vu All Over Again?}
With so many proofs of the central theorem already on the books, similarities are bound to come up. A reader of an earlier version of this paper found similarity with \cite{gilbert-1963}. The similarity, however, is confined to setup stage of the proof. To make an analogy with the game of chess, the popular \emph{Ruy Lopez} opening has variations extending to move nine; there are substantially different chess games whose first nine moves are identical. The setup portions of our proof prepare for application of the Law of sines, which is covered in most high school geometry courses. The proof in \cite{gilbert-1963} relies on the concept of  \emph{concyclic quadrilaterals}, not covered in that curriculum. See the books \cite{hartshorne-2000} and \cite{lee-2014}, in which the Law of sines is mentioned but \emph{concyclic} does not appear in the index. M.~Lewin points out \cite{lewin-1974} that the proof in \cite{gilbert-1963}, though selected from manifold contributions, is ``essentially identical'' to C.L.~Lehmus' original 1850 proof. The historically oriented and much-valued book \emph{Geometry Through Its History} \cite{ostermann-2012} gives a proof, attributed to V. Th\'ebault (see \cite{thebault-1937}), which is basically the same as \cite{gilbert-1963}, though with more details provided.   Steiner's 1840s proof is also provided in \cite{ostermann-2012}, but we could not find mention of Lehmus' 1850 proof. It's quite likely that the topological space of Steiner-Lehmus proofs is a loop space: ``The circle closed'' \cite{lewin-1974} .

\section{Why Reprove it?}
Some theorems engender a steady stream of new, or newly presetented proofs.  But why?  The Morley Trisection Theorem (MTT) \cite{oakley-1978} is a case in point. It shares with the SLT the feature of being a triangle theorem whose statement is at the level  of ancient geometry; both theorems are included in Martin Gardner's compendium \cite{gardner-1961}. Proofs  of the MTT abound; see \cite{oakley-1978} for a multitude (plethora?), and Math Reviews  for dozens more which followed. Yet the motivation for reproof is subtle. Perhaps a feeling of insecurity is at play. After a lost puppy is found and returned, the owners may check and recheck to ensure that the puppy is still there and not lost again. The MTT was discovered accidentally; it is a surprising result, and has been called a miracle. Indeed, the MTT surprise element does not seem to fade, even after repeated washings. Re-proofs may serve to reassure that the theorem is still there, valid and not forgotten. (Of course, new proofs add value in other ways as well, e.g,  by providing opportunities for deeper and broader understanding.)

For the SLT the standard explanation for reproofs  is that a \emph{direct} argument is desired; if one cannot be found then the goal is to approach directness ever more closely.  J.H.~Conway and A.~Ryba have argued eloquently that a direct proof is not possible \cite{conway-2014}, at least not in certain familiar contexts of geometry. More recently. V.~Pambuccian proved the abstract existence of a direct proof \cite{pambuccian-2017}, though one of perhaps several thousand lines in length, and  A.~Kellison gave a concrete, \emph{machine checked}, direct proof \cite{kellison-2021}. These proofs invoke adapted, context-shifting axioms; they do not invalidate the Conway \& Ryba dictum. And we are sure to see more proofs in the future. The experience of the four 
\texttt{\footnotesize [make that five]} 
aforementioned prominent mathematicians suggests that, in addition to directness of proof, interest in the SLT is sustained by another reason: it is a career launching pad for budding mathematicians.

\noindent
{UMass Boston and   University of New Hampshire} \\ \\
\url{eric.grinberg@umb.edu} \, and \, \url{mo@unh.edu}

\end{document}